\documentclass[11pt]{amsart}
\usepackage[mathscr]{euscript}
\usepackage{amssymb}

\usepackage{amscd}
\usepackage{amsmath, amssymb}
\usepackage{amsfonts}
\usepackage[colorlinks,linkcolor=blue,citecolor=blue, pdfstartview=FitH]{hyperref}
\usepackage[all]{xy}

\numberwithin{equation}{section}

\theoremstyle{plain}
\newtheorem{thm}{Theorem}[section]
\newtheorem{theorem}[thm]{Theorem}
\newtheorem{lemma}[thm]{Lemma}
\newtheorem{corollary}[thm]{Corollary}
\newtheorem{proposition}[thm]{Proposition}

\theoremstyle{definition}

\newtheorem{definition}[thm]{Definition}

\newtheorem{defn-thm}[thm]{Definition-Theorem}

\newcommand{\im}{{ \mathrm{im}\,}}

\newcommand{\C}{{ \mathbb{C} }}
\newcommand{\Z}{{ \mathbb{Z} }}

\newcommand{\G}{{\textrm{Grass}\,}}
\newcommand{\p}{{ \partial }}
\newcommand{\pb}{{ \bar{\partial} }}

\newcommand{\wg}{{\,\wedge\,}}
\newcommand{\Hom}{{ \textrm{Hom}\,}}

\begin{document}

\title{Variation of Hodge structures for non-K\"ahler manifolds}
\author{Wei Xia}
\address{Wei Xia, Mathematical Science Research Center, Chongqing University of Technology, Chongqing, P.R.China, 400054.} \email{xiawei@cqut.edu.cn, xiaweiwei3@126.com}

\thanks{This work was supported by the National Natural Science Foundation of China No. 11901590 and Natural Science Foundation of Chongqing (China) No. CSTB2022NSCQ-MSX0876. This work is also partially supported by the Natural Science Foundation of Chongqing (China) No. CSTB2024NSCQ-LZX0040 and No. CSTB2023NSCQ-LZX0042.}
\date{\today}

\begin{abstract}
In this paper, we discuss unpolarized, complex variation of Hodge structures for non-K\"ahler manifolds. In particular, consider a holomorphic family of compact complex manifolds with central fiber $X$. Suppose that the inclusion of complexes $(F^pA^\bullet(X),d)\hookrightarrow (A^\bullet(X),d)$ induces injective maps
\[
H^j(F^pA^\bullet(X),d)\longrightarrow H^j(A^\bullet(X),d),\qquad j=k,\ k+1.
\]
Under this condition, we show that the period map is holomorphic and transversal. An essential ingredient for the proof is a deformation theory of filtered $(p,q)$-forms. We prove a jumping formula for filtered de Rham cohomology and obtain a necessary and sufficient condition for its dimension to be invariant under deformation of complex structures.
\end{abstract}

\maketitle

%\tableofcontents
%%%%%%%%%%%%%%%%%%%%%%%%%%%%%%%%%%%%%%%%%%%%%%%%%%%%%%%%%%%
\section{Introduction}
The theory of variation of Hodge structures (VHS for short) is an extensively studied object in complex geometry and related subjects (see \cite{GGK13,CMP17,KP16,CK99} and the references therein). The purpose of this paper is to extend some parts of Griffiths' classical theory of VHS \cite{Gri68} to the non-K\"ahler setting, see \cite{Sim91,Kir15,Kas21} for extensions in other directions. There are two motivations to do this, the first is the works of Popovici and Anthes-Cattaneo-Rollenske-Tomassini \cite{Pop19H,ACRT18}, where they showed local Torelli theorem holds for Calabi-Yau $\p\pb$-manifolds and $\p\pb$-complex symplectic manifolds, respectively; the second is that we want to understand Voisin's proof (see \cite{Voi05}) of a density criterion for complex projective manifolds\footnote{This has been generalized by Rutong Chen \cite{Chen22} to Moishezon manifolds by using the results obtained in this paper.}.

Let $X$ be a compact complex manifold. There is the \emph{Hodge filtration on forms}
\[
F^pA^{p+q}(X)=\bigoplus_{i=0}^{\infty}A^{p+i,q-i}(X),\quad p,q\in\Z,
\]
where $A^{p,q}(X)$ is the space of smooth complex valued differential forms on $X$. Following \cite{WX23}, elements in $F^pA^{p+q}(X)$ will be called \textit{filtered $(p,q)$-forms on $X$}. Since Hodge filtrations are compatible with the exterior differential $d$, it induces \emph{Hodge filtrations on the de Rham cohomology}:
\[
F^pH^{p+q}(X):=\im \left( \frac{\ker d\cap F^pA^{p+q}(X) }{\im d\cap F^pA^{p+q}(X) }\longrightarrow H^{p+q}(X,\C) \right),\quad p,q\in\Z.
\]
As in the classical theory, we are interested in the variations of $F^pH^{p+q}(X)$ when the complex structure on $X$ deforms. For this purpose, we first establish a deformation theory for filtered $(p,q)$-forms.

More precisely, let
\[
\pi: (\mathcal{X}, X)\to (B,0)
\]
be a complex analytic family over a small polydisc $B\subset\C^m$ such that for each $t\in B$ the complex structure on $X_t$ is represented by a Beltrami differential $\phi(t)$. Given
\[
\alpha_0\in\ker d\cap F^{p}A^{p+q}(X),
\]
and $T\subseteq B$, which is an analytic subset of $B$ containing $0$, a \textit{deformation} of $\alpha_0$ (w.r.t. $\pi$ ) on $T$ is a family of forms
\[
\alpha (t)\in F^{p}A^{p+q}(X),
\]
such that
\begin{itemize}
  \item[(1)] $\alpha (t)$ is holomorphic in $t$ and $\alpha(0) = \alpha_0$;
  \item[(2)] for any $t\in T$, we have
  \begin{equation}
  d_{\phi(t)}\alpha (t)=0.
  \end{equation}
\end{itemize}
Here $d_{\phi(t)}=d-\mathcal{L}^{1,0}_{\phi(t)}$, where $\mathcal{L}^{1,0}_{\phi(t)}=i_{\phi(t)}\partial-\partial i_{\phi(t)}$. Locally, write $\phi(t)=\sum_a\eta_a\otimes V_a$ with $(0,1)$-forms $\eta_a$ and $(1,0)$-vector fields $V_a$. Then $i_{\phi(t)}\alpha=\sum_a\eta_a\wedge\iota_{V_a}\alpha$, where $\iota_{V_a}$ is contraction by $V_a$. See \eqref{eq-LRY15} for the conjugation formula.

For any given closed filtered $(p,q)$-form $\alpha_0\in\ker d\cap F^{p}A^{p+q}(X)$, the deformations $\alpha(t)$ of $\alpha_0$ are by no means unique. Motivated by the deformation theory of Dolbeault cohomology classes \cite{Xia22}, we are led to consider canonical deformations, i.e. we will construct deformations $\alpha(t)$ by using Hodge decomposition for filtered forms with respect to a fixed Hermitian metric on $X$.

Assume now that a fixed Hermitian metric has been equipped on $X$. A version of Hodge decomposition holds for filtered $(p,q)$-forms (cf. Subsection \ref{subsec-Hodge-decomp}) with the formal adjoint of $d$, Laplacian and Green operator being $d^*_p, \triangle_{p}, G_p$, respectively. A deformation $\alpha(t)$ of $\alpha_0$ is \emph{canonical} if it is $d^*_p$-closed up to a $d_{\phi(t)}$-exact deformation (see Definition \ref{def-deformation-Dol-Zrpq}).
We say that the \textit{deformations of filtered $(p,q)$-forms on $X$ are (canonically) unobstructed} (w.r.t. $\pi$) if for any $\alpha_0\in\ker d\cap F^{p}A^{p+q}(X)$
there is a (canonical) deformation of $\alpha_0$ on $B$. Clearly, canonical unobstructedness implies unobstructedness.

It may happen that $\alpha_0$ is not $d$-exact while $\alpha(t)$ is $d_{\phi(t)}$-exact for some $t$, i.e. $\alpha(t)$ is not a faithful deformation in the sense of Definition \ref{def-deformation-Dol-Zrpq}. As in the case of deformations of Dolbeault cohomology, this together with possible obstructions are exactly the causes which make $\dim F^pH^{p+q}(X_t)$ jump at $t=0$. In Theorem \ref{thm-jumping-formula}, we will prove the following: for any $0\leq p\leq k$ and $t\in B$, set
\[
v^{p,k}_t:=\dim \ker\triangle_{p}\cap F^{p}A^{k}(X)-\dim F^pA^{k}(X)\cap\ker d_{\phi(t)}\cap\ker d^*_p \geq 0,
\]
where $d^*_p=d^*-\p^*\Pi^{p,\bullet-p}$ is the formal adjoint of $d$ on $F^{p}A^{k}(X)$ and $\Pi^{p,\bullet-p}$ projects a form of degree $k$ onto its $(p,k-p)$ component (see Subsection \ref{subsec-Hodge-decomp}). Then we have $(v^{p,-1}_t:=0)$
\begin{equation}\label{eq-jumping-formula-0}
\dim F^pH^{k}(X)=\dim F^pH^{k}(X_t)+v^{p,k}_t+v^{p,k-1}_t.
\end{equation}
This may be called \emph{the jumping formula for filtered de Rham cohomology}, see \cite{Xia22,HX25} for versions of Dolbeault/Bott-Chern/Aeppli cohomology. This formula may be viewed as a refinement of the classical result obtained by Kodaira-Spencer (see formula $(43)$ in \cite[pp.\,68]{KS60}). In their formula, $v^q_t$ is defined via the spectrum of Laplacian operators. Here, the integers $v^{p,k}_t$ have significant geometric meanings, roughly speaking, it measures the size of classes in $F^pH^{k}(X)$ whose canonical deformation does not exist on $t$. As direct applications of the jumping formula \eqref{eq-jumping-formula-0}, we get the following
\begin{corollary}
The following holds:
\begin{enumerate}
  \item $\dim F^pH^{k}(X_t)$ is independent of $t\in B$ if and only if the deformations of filtered $(p,k-p)$-forms and filtered $(p,k-1-p)$-forms are canonically unobstructed;
  \item For any $0\leq p\leq k$, the alternating sum $\sum_{i=0}^k (-1)^{k-i}\dim F^pH^{i}(X_t)$ and $\dim E^{0,q}_1(X_t)/E^{0,q}_\infty(X_t)$ are upper semi-continuous functions of $t\in B$ (in analytic Zariski topology). Moreover, the deformations of filtered $(p,k-p)$-forms are canonically unobstructed if and only if $\sum_{i=0}^k (-1)^{k-i}\dim F^pH^{i}(X_t)$ is independent of t;
  \item The notion of canonically unobstructedness for filtered forms are independent of the choices of Hermitian metrics.
\end{enumerate}
\end{corollary}
In $(2)$, $E^{p,q}_r(X_t)$ are the spaces in the Fr\"olicher spectral sequence of $X_t$. The upper semi-continuity of $\dim E^{0,q}_1(X_t)/E^{0,q}_\infty(X_t)$ for $q=1$ or $n-1$ was shown by Stelzig \cite{Ste22} (in Euclidean topology). On the other hand, it was shown by Flenner \cite{Fle81} that for any $q$, $\sum_{i=0}^{q}(-1)^{q-i}\dim H^{i}(X_t,\mathcal{E}\mid_{X_t})$ is upper semi-continuous in analytic Zariski topology, where $\mathcal{E}$ is a coherent analytic sheaf on the total space $\mathcal{X}$, see also \cite[Thm.\,5.10]{BDIP02}. Notice also that when $p=0, k=2\dim_\C X$, the alternating sum $\sum_{i=0}^k \dim (-1)^{k-i}F^pH^{i}(X_t)=\sum_{i=0}^k \dim (-1)^{i}H^{i}(X_t)$ is just the topological Euler number of $X$.

Let $\pi: (\mathcal{X}, X)\to (B,0)$ be a complex analytic family over a small polydisc $B\subset\C^m$ such that for each $t\in B$ the complex structure on $X_t$ is represented by a Beltrami differential $\phi(t)$. Following \cite[Chp.\,10]{Voi02I}, the \emph{period map} is defined as follows
\begin{equation}\label{eq-period-map-0}
\Phi^{p,k}:B\longrightarrow \G (f^{p,k}, H^{k}(X,\C)),\quad t\longmapsto F^pH^{k}(X_t),
\end{equation}
where $f^{p,k}:=\dim F^pH^{k}(X)$. Our main result is the following:
\begin{theorem}[=Theorem \ref{thm-period-map}]\label{thm-period-map-0}
Fix $p\leq k$ and assume
\begin{equation}\label{eq-assumption-1}
  F^{p}A^{k+1}(X)\cap dA^{k}(X) =dF^{p}A^{k}(X),
\end{equation}
and
\begin{equation}\label{eq-assumption-2}
  F^{p}A^{k}(X)\cap dA^{k-1}(X) =dF^{p}A^{k-1}(X).
\end{equation}
The following holds:
\begin{enumerate}
  \item The period map $\Phi^{p,k}$ is holomorphic;
  \item Griffiths transversality: the tangent map
  \[
  d\Phi^{p,k}_0:T_0B\longrightarrow  \Hom(F^pH^{k}(X), H^{k}(X,\C)/F^pH^{k}(X))
  \]
  has values in $\Hom(F^pH^{k}(X), F^{p-1}H^{k}(X)/F^pH^{k}(X))$;
\end{enumerate}
Let now $\omega$ be a Hermitian metric on $X$. Denote by $\mathcal H_p$ the orthogonal projection onto the space of harmonic filtered forms $\ker\triangle_p$ in $F^pA^\bullet(X)$ (see Subsection \ref{subsec-Hodge-decomp}). If $\mathcal{H}_p=\mathcal{H}_{p-1}\mid_{F^{p}A^{\bullet}(X)}$ and $\dim B=1$, then the tangent map of $\Phi^{p,k}$ is given by
\[
d\Phi^{p,k}_0(\frac{\p}{\p t}\Big|_{t=0})([\alpha_0])=[\mathcal{H}_{p-1}i_{\phi_1}\alpha_0]\in F^{p-1}H^{k}(X),\quad [\alpha_0]\in F^pH^{k}(X),
\]
where $t$ is a local coordinate function on $B$ and $\phi_1$ is the first order coefficient in $\phi(t)=\sum_{j\geq 1}\phi_jt^j$.
\end{theorem}
The assumption \eqref{eq-assumption-1} says exactly that the inclusion of complexes
\[
(F^pA^\bullet(X),d)\hookrightarrow (A^\bullet(X),d)
\]
induces an injective map
\[
H^{k+1}(F^pA^\bullet(X),d)\longrightarrow H^{k+1}(A^\bullet(X),d).
\]
Set $q=k-p$ for the spectral-sequence indices. In terms of the differentials $d_r^{p,q}:E_r^{p,q}(X)\longrightarrow E_r^{p+r,q-r+1}(X)$, \eqref{eq-assumption-1} is equivalent to $\bigoplus_{r\geq i\geq 1}d_r^{p-i,q+i}=0$, see Corollary B.6 in \cite{WX23}. For example, for the Iwasawa manifold, \eqref{eq-assumption-1} and \eqref{eq-assumption-2} (being equivalent to $d_1^{p-1,q+1}=0$ and $d_1^{p-1,q}=0$, respectively) hold whenever $p\neq2$. For more examples of this type, see \cite[Sec.\,5.2]{WX23} and \cite[Sec.\,4]{COUV16}.

A main ingredient to prove this theorem is the fact that the exponential operator preserves Hodge filtrations, see Theorem \ref{thm-exp-FpH}. The period map defined by \eqref{eq-period-map-0} concerns unpolarized, complex VHS for non-K\"ahler manifolds. It is an interesting problem to extend $\Phi^{p,k}$ to the polarized version. Unfortunately, in our case as K\"ahler metrics are not assumed to exist, the Hodge-Riemann bilinear relations and hard Lefschetz theorem are not available. Moreover, it is perhaps worthwhile to point out that if the Hermitian metric $h$ is K\"ahler, the condition $\mathcal{H}_p=\mathcal{H}_{p-1}\mid_{F^{p}A^{\bullet}(X)}$ automatically holds for all $p\geq0$ (see Proposition \ref{prop-kahler}). It seems an interesting question to determine the type of special Hermitian metrics (see \cite{Fin24,Popbook} and the references therein for detailed discussions) that make $\mathcal{H}_p=\mathcal{H}_{p-1}\mid_{F^{p}A^{\bullet}(X)}$ hold for some fixed $p\geq0$.

\vskip 1\baselineskip \textbf{Statements and Declarations.}
The author has no competing interests to declare that are relevant to the content of this article.

\vskip 1\baselineskip \textbf{Acknowledgements.} We would like to thank Prof. Kefeng Liu for useful comments. Many thanks to Xueyuan Wan, Rutong Chen, Yang Shen and Ya Deng for stimulating discussions. This work was supported by the National Natural Science Foundation of China No. 11901590 and Natural Science Foundation of Chongqing (China) No. CSTB2022NSCQ-MSX0876. This work is also partially supported by the Natural Science Foundation of Chongqing (China) No. CSTB2024NSCQ-LZX0040 and No. CSTB2023NSCQ-LZX0042.

\section{The exponential operator} \label{}
The exponential operator\footnote{An implicit form of the exponential operator was used by Griffiths in \cite[Prop.\,1.11]{Gri68}.}
\[
e^{i_{\phi(t)}}:=\sum_{k=0}^{\infty} \frac{i_{\phi(t)}^k}{k!}
\]
which was introduced in the work of Todorov (see Remark $2.4.9.$ in ~\cite[p.\,339]{Tod89}) turns out to be very useful for studying deformations of complex manifolds, see e.g. \cite{Cle05,LSY09,LRY15,RZ18,LS20,Xia19dBC}.
Let $\pi: (\mathcal{X}, X)\to (B,0)$ be a complex analytic family in the sense of Kodaira-Spencer, that is, $\pi$ is a holomorphic submersion which is surjective and proper over the connected complex manifold $B$. Assume $B$ is a small polydisc (centered in $0$) with coordinates $t=(t_1,t_2,\cdots)$, then by Ehresmann's theorem we have the following commutative diagram
\[
  \xymatrix{
  X\times B \ar[dr] \ar[r]^-{F}
                & \mathcal{X} \ar[d]^-{\pi}  \\
                & B~, }
\]
where $F$ is a diffeomorphism. For any $t\in B$ set $f_t:=F\mid_{X\times \{t\}}$, then $f_t$ is a diffeomorphism from $X$ to $X_t$. We use the isomorphism
\[
f_t^*:A^\bullet(X_t)\longrightarrow A^\bullet(X)
\]
to identify forms/cohomologies on $X_t$ with forms/cohomologies on $X$.
Let $z^1,\cdots, z^n$ and $w^1,\cdots, w^n$ be holomorphic coordinates on $X$ and $X_t$ respectively, the Beltrami differential can be defined by
\begin{equation}\label{Beltrami differential coord}
\phi(t):=\left(\frac{\partial w}{\partial z}\right)^{-1}_{\alpha\gamma}\frac{\partial w^{\gamma}}{\partial \bar{z}^{\beta}} d\bar{z}^{\beta}\otimes\frac{\partial}{\partial z^{\alpha}}\in A^{0,1}(X, T^{1,0}X)~,
\end{equation}
where by abuse of notations, we write $w^i=f_t^*w^i=w^i\circ f_t$ for each $i=1,2,\cdots,n$.
We note that
\begin{equation}\label{dw in terms of dz}
d w^{\beta}=\frac{\partial w^{\beta}}{\partial z^{\alpha}}dz^{\alpha}+\frac{\partial w^{\beta}}{\partial \bar{z}^{\alpha}}d\bar{z}^{\alpha}
=\frac{\partial w^{\beta}}{\partial z^{\alpha}}(1+i_{\phi(t)})dz^{\alpha},
\end{equation}
and the fact that $e^{i_{\phi(t)}}$ is compatible with wedge product (cf. \cite[p.\,339]{Cle05}) which implies the following
\begin{align*}
e^{i_{\phi(t)}} (dz^{i_1}\wg\cdots\wg dz^{i_p})
&=(e^{i_{\phi(t)}}dz^{i_1}) \wg \cdots \wg (e^{i_{\phi(t)}}dz^{i_p} )\\
&=(dz^{i_1}+i_{\phi(t)}dz^{i_1}) \wg \cdots \wg (dz^{i_p}+i_{\phi(t)}dz^{i_p}),
\end{align*}
where for each $k=1,\cdots,p$,
\[
i_{\phi(t)}^2dz^{i_k}=0\Longrightarrow
e^{i_{\phi(t)}}dz^{i_k}
=dz^{i_k}+i_{\phi(t)}dz^{i_k}.
\]
As a result, we see that
\begin{equation}\label{eq-exp}
e^{i_{\phi(t)}}~: A^{p,0}(X)\longrightarrow A^{p,0}(X_t).
\end{equation}
We will use the following formula (see e.g. \cite{LRY15,FM06,Xia19deri} and \cite[p.\,334]{Ma05}):
\begin{equation}\label{eq-LRY15}
d_{\phi(t)}:=e^{-i_{\phi(t)}}de^{i_{\phi(t)}} = \partial + \bar{\partial}_{\phi(t)}=\partial + \bar{\partial}-\mathcal{L}_{\phi(t)}^{1,0}~,
\end{equation}
where $\mathcal{L}_{\phi(t)}^{1,0}:=i_{\phi(t)}\p-\p i_{\phi(t)}$ is the Lie derivative.
A remarkable property of the exponential operator (as observed in \cite{FM06,FM09,WZ20}) is that while it does not preserve all the types of $(p,q)$-forms, it does map elements in a filtration on $X$ to elements in the corresponding filtration on $X_t$:
\begin{equation}\label{eq-weidingchang}
e^{i_{\phi(t)}}~: F^pA^{k}(X)=\oplus_{p\leq\lambda\leq n} A^{\lambda,k-\lambda}(X)\longrightarrow F^pA^{k}(X_t),
\end{equation}
where $n=\dim_\C X$. In fact, let $\varphi\in A^{\lambda,k-\lambda}(X)$ with $p\leq\lambda\leq k$, locally we write $\varphi=\varphi_{IJ}dz^I\wedge d\bar{z}^J$ where $dz^I = dz^{i_1}\wedge\cdots\wedge dz^{i_\lambda}$ and $d\bar{z}^J = d\bar{z}^{j_1} \wedge\cdots\wedge d\bar{z}^{j_{k-\lambda}}$, then
\begin{align*}
&e^{i_{\phi(t)}}\varphi\\
=&\varphi_{IJ}(e^{i_{\phi(t)}}dz^I)\wedge d\bar{z}^J\\
=&\varphi_{IJ}(e^{i_{\phi(t)}}dz^I)\wg (\frac{\p\bar{z}^{j_1}}{\p w^l}dw^l+\frac{\p\bar{z}^{j_1}}{\p \bar{w}^l}d\bar{w}^l)\wg\cdots\wg (\frac{\p\bar{z}^{j_{k-\lambda}}}{\p w^l}dw^l+\frac{\p\bar{z}^{j_{k-\lambda}}}{\p \bar{w}^l}d\bar{w}^l).
\end{align*}
Notice that the factor $e^{i_{\phi(t)}}dz^I$ has type $(\lambda,0)$ on $X_t$ by \eqref{eq-exp}. Thus every term in the wedge product has holomorphic degree at least $\lambda\geq p$ on $X_t$, even if every remaining factor contributes only a $d\bar w^l$ term. Hence $e^{i_{\phi(t)}}\varphi$ belongs to $F^pA^k(X_t)$. It follows that (cf. \cite{WX23})
\begin{equation}\label{eq-wdc-dc}
e^{i_{\phi(t)}}:(F^\bullet A^{\bullet}(X),d_{\phi(t)} )\to (F^\bullet A^{\bullet}(X_t),d)
\end{equation}
is an isomorphism of filtered complexes.

\section{Deformations of filtered forms} \label{Deformations of filtered froms}
In this section, we establish a deformation theory for filtered forms which is modeled on \cite{Xia22}.
\subsection{Hodge decomposition for filtered forms}\label{subsec-Hodge-decomp}
Assume $X$ has been equipped with a fixed Hermitian metric. There is a Hodge decomposition for filtered forms (cf. Appendix A of \cite{WX23}). In fact, set
\[
d^*_p=\Pi^{\geq p}d^*:F^{p}A^{\bullet}(X)\longrightarrow F^{p}A^{\bullet}(X),
\]
where $\Pi^{\geq p}:A^{\bullet}(X)\to F^pA^{\bullet}(X)$ is the linear projection onto $F^pA^{\bullet}(X)$. It turns out $d^*_p$ is the formal adjoint of $d$ for filtered forms, namely,
\[
(d\alpha,\beta)=(\alpha,d^*_p\beta),\quad \text{for~any}~\alpha\in F^pA^{\bullet}(X), \beta\in F^pA^{\bullet+1}(X).
\]
Notice that $d^*_p=d^*-\p^*\Pi^{p,\bullet-p}$ where $\Pi^{p,q}:A^{\bullet}(X)\to A^{p,q}(X)$ is the linear projection onto $A^{p,q}(X)$. So we have $d^*_p\mid_{F^{p+1}A^{\bullet}(X)}=d^*$.
The \emph{Laplacian operator for filtered forms} is defined by
\[
\triangle_p=dd^*_p+d^*_pd:F^{p}A^{\bullet}(X)\longrightarrow F^{p}A^{\bullet}(X),
\]
then there is the \emph{Green operator} of $\triangle_p$, such that
\[
1=\mathcal{H}_p+\triangle_pG_p=\mathcal{H}_p+G_p\triangle_p,\quad\text{on}~F^pA^{\bullet}(X),
\]
where $\mathcal{H}_p$ (we also write $\mathcal{H}_p^k$ for definiteness) is the projection onto the harmonic space $\ker\triangle_{p}\cap F^{p}A^{k}(X)$.
Now let $\pi: (\mathcal{X}, X)\to (B,0)$ be a complex analytic family over a small polydisc $B\subset\C^m$ such that for each $t\in B$ the complex structure on $X_t$ is represented by a Beltrami differential $\phi(t)$. We have the following:
\begin{itemize}
  \item For any $t\in B$, set
  \[
  \triangle_{p\phi(t)}=d_{p\phi(t)}^*d_{\phi(t)}+d_{\phi(t)} d_{p\phi(t)}^*:F^{p}A^{\bullet}(X)\longrightarrow F^{p}A^{\bullet}(X),
  \]
  where $d_{p\phi(t)}^*=\Pi^{\geq p}d_{\phi(t)}^*$ and $d_{\phi(t)}^*$ is the formal adjoint of $d_{\phi(t)}$;
  \item There is a Green operator $G_{p\phi(t)}:F^{p}A^{\bullet}(X)\longrightarrow F^{p}A^{\bullet}(X)$ such that
  \[
  1=\mathcal{H}_{p\phi(t)}+\triangle_{p\phi(t)}G_{p\phi(t)}
   =\mathcal{H}_{p\phi(t)}+G_{p\phi(t)}\triangle_{p\phi(t)},
  \]
  where $\mathcal{H}_{p\phi(t)}$ is the harmonic projection operator.
\end{itemize}
It is clear that $\triangle_{p\phi(t)}$ is a differentiable family of formally self-adjoint, elliptic differential operators in the sense of Kodaira-Spencer \cite{KS60} (see also \cite{Kod86}). As a result, by Kodaira-Spencer's upper semi-continuity theorem, we get that for any $t\in B$,
\[
\dim F^pH^{k}(X)=\dim \ker \triangle_{p}\geq \dim \ker \triangle_{p\phi(t)}=\dim F^pH^{k}_{d_{\phi(t)}}(X)= \dim F^pH^{k}(X_t),
\]
where $F^pH^{k}_{d_{\phi(t)}}(X):=\frac{\ker d_{\phi(t)}\cap F^pA^{k}(X)}{\im d_{\phi(t)}\cap F^pA^{k}(X)}$ and it follows from \eqref{eq-weidingchang} that
\[
F^pH^{k}_{d_{\phi(t)}}(X)\cong F^pH^{k}(X_t).
\]
As a result, we have the following:
\begin{proposition}\label{prop-upper-semicontinuous}
Let $\pi: (\mathcal{X}, X)\to (B,0)$ be a complex analytic family over a small polydisc $B\subset\C^m$. Then $\dim F^pH^{\bullet}(X_t)$ is an upper semi-continuous function of $t\in B$.
\end{proposition}
An immediate consequence is that for any $q\geq0$, $\dim E^{0,q}_\infty$ is lower semi-continuous in deformation of complex structures. In fact, from
\begin{equation}\label{eq-E0q}
E^{0,q}_\infty(X)\cong \frac{F^0H^{q}(A^{\bullet}(X),d)}{F^{1}H^{q}(A^{\bullet}(X),d)}\cong \frac{H^{q}_{dR}(X)}{F^{1}H^{q}(X) },
\end{equation}
we get $\dim E^{0,q}_\infty(X)=b^q-\dim F^{1}H^{q}(X)$, where $b^q=\dim H^{q}_{dR}(X)$ is the $q$-th Betti number of $X$. By using Proposition \ref{prop-upper-semicontinuous}, we get the following
\begin{corollary}
$\dim E^{0,q}_\infty$ is lower semi-continuous in deformation of complex structures. $\dim E^{0,q}_1/E^{0,q}_\infty$ is upper semi-continuous in deformation of complex structures.
\end{corollary}
\begin{proof}The second statement follows since $\dim E^{0,q}_1=h^{0,q}$, the $(0,q)$-th Hodge number, is upper semi-continuous and the sum of two upper semi-continuous functions is still upper semi-continuous.
\end{proof}
There is a strengthened version of this result which holds in analytic Zariski topology, see Corollary \ref{Coro-123}.
\subsection{Construction of deformations}
We write $\phi(t)=\sum_{i\geq1}\phi_i$.
\begin{proposition}\label{prop-unique-power-seris-solution}
\begin{enumerate}
  \item For any $t\in B$, the following natural homomorphism induced by the inclusion $\ker d_{\phi(t)}\cap\ker d^*_p\subset\ker d_{\phi(t)}$ is an isomorphism:
\[
\frac{F^pA^{k}(X)\cap\ker d_{\phi(t)}\cap\ker d^*_p}{F^pA^{k}(X)\cap\im d_{\phi(t)}\cap\ker d^*_p}\longrightarrow
\frac{F^pA^{k}(X)\cap\ker d_{\phi(t)}}{F^pA^{k}(X)\cap\im d_{\phi(t)}}=F^pH^{k}_{d_{\phi(t)}}(X).
\]
  \item For any given $\alpha\in F^{p}A^{k}(X)$, if $d_{\phi(t)}\alpha = d\alpha - \mathcal{L}_{\phi(t)}^{1,0}\alpha =0$ and $d^*_p\alpha=0$, then we must have
\[
\alpha = \mathcal{H}_{p}\alpha + d^*_p G_p\mathcal{L}_{\phi(t)}^{1,0}\alpha .
\]
  \item For any fixed $\alpha_0\in\ker\triangle_{p}\cap F^{p}A^{k}(X)$, the equation
\begin{equation}\label{Kuranishi eq}
\alpha=\alpha_0 + d^*_p G_p\mathcal{L}_{\phi(t)}^{1,0}\alpha,
\end{equation}
has an unique solution given by $\alpha=\alpha(t)=\sum_{k\geq0} \alpha_k$ and
\begin{equation}\label{eq-canonical-powerseries}
\alpha_k=d^*_p G_p\sum_{i+j=k} \mathcal{L}_{\phi_i}^{1,0}\alpha_{j}\in F^{p}A^{k}(X),\quad k>0,
\end{equation}
which converges for $|t|$ small.
\item Let $\alpha$ be a solution of the equation \eqref{Kuranishi eq}.
Then for any $t\in B$, we have
\begin{equation}\label{eq-equivalence-harmonic part vanish}
d_{\phi(t)}\alpha=0 \Leftrightarrow \mathcal{H}_p\mathcal{L}_{\phi(t)}^{1,0}\alpha =0.
\end{equation}
\item Set
\begin{align*}
\hat{g}_t:\ker\triangle_{p}\cap F^{p}A^{k}(X)&\longrightarrow \ker d_p^*\cap\im d_{\phi(t)}\cap F^{p}A^{k+1}(X)\\
x_0&\longmapsto d_{\phi(t)}x(t),
\end{align*}
where $x(t)$ is the unique solution of $x(t)=x_0+d^*_p G_p\mathcal{L}_{\phi(t)}^{1,0}x(t)$. Then the following holds:
\[
\dim\ker\triangle_{p}\cap F^{p}A^{k}(X)=\dim\left(\ker d_p^*\cap\ker d_{\phi(t)}\right)^k+
\dim \left(\ker d_p^*\cap\im d_{\phi(t)}\right)^{k+1},
\]
where
\begin{itemize}
  \item $\left(\ker d_p^*\cap\ker d_{\phi(t)}\right)^k:=\ker d_p^*\cap\ker d_{\phi(t)}\cap F^{p}A^{k}(X)$;
  \item $\left(\ker d_p^*\cap\im d_{\phi(t)}\right)^{k+1}:=\ker d_p^*\cap\im d_{\phi(t)}\cap F^{p}A^{k+1}(X)$.
\end{itemize}
\end{enumerate}
\end{proposition}
\begin{proof}
This follows from minor modifications of the arguments in \cite[Sec.\,4]{Xia22}. Indeed, by replacing $F^pA^{k}(X), d_{\phi(t)}, d_p^*, \triangle_{p}, G_p, \mathcal{H}_{p}$ with the spaces and operators appeared above: $A^{0,q}(X, E), \pb_{\phi(t)}, \pb^*, \Box, G_{\phi(t)}, \mathcal{H}_{\phi(t)}$, correspondingly, the proof of $(1)-(5)$ are essentially contained in Proposition 4.5 (for $(1)$), Proposition 4.1 (for $(2),(3)$), Proposition 4.2 (for $(4)$), Proposition 4.6 (for $(5)$) of \cite{Xia22}, respectively.
\end{proof}
\begin{proposition}\label{prop-obstruction-form}
For any $\alpha_0\in\ker\triangle_{p}\cap F^{p}A^{k}(X)$, we have $\mathcal{L}_{\phi(t)}^{1,0}\alpha(t)\in\ker d$ where $\alpha(t)=\sum_{m\geq0} \alpha_m$ is given by \eqref{eq-canonical-powerseries}.
\end{proposition}
\begin{proof} By applying Hodge decomposition theorem to $\mathcal{L}_{\phi(t)}^{1,0}\alpha(t)$ we get
\[
-d_{\phi(t)}\alpha(t)=-(d-\mathcal{L}_{\phi(t)}^{1,0})\alpha(t)
=\mathcal{H}_p\mathcal{L}_{\phi(t)}^{1,0}\alpha(t)+
d\big(d^*_pG_p\mathcal{L}_{\phi(t)}^{1,0}\alpha(t)-\alpha(t)\big)
+d^*_pdG_p\mathcal{L}_{\phi(t)}^{1,0}\alpha(t).
\]
It follows from $(4)$ of Proposition \ref{prop-unique-power-seris-solution} that $d_{\phi(t)}\alpha(t)=0\Leftrightarrow\mathcal{H}_p\mathcal{L}_{\phi(t)}^{1,0}\alpha(t)=0$ which implies
\[
d^*_pdG_p\mathcal{L}_{\phi(t)}^{1,0}\alpha(t)=0\Longrightarrow
\mathcal{L}_{\phi(t)}^{1,0}\alpha(t)\in\ker d.
\]
\end{proof}
In view of Proposition \ref{prop-obstruction-form}, $[\mathcal{L}_{\phi(t)}^{1,0}\alpha(t)]\in F^pH^{k+1}(X)$ may be called the \emph{obstruction class} of $\alpha_0$.
\begin{definition}\label{def-V_t^k}
For any $t\in B$ and a vector subspace $V\subseteq \ker\triangle_{p}\cap F^{p}A^{k}(X)$, we set
\begin{align*}
V_{t}^{p,k}:=
\Big\{ \alpha_0\in V \mid~~ &d_{\phi(t)}\alpha(t)=0~\text{where}
~\alpha(t)\\
&\text{is the unique solution of}~\alpha (t)=\alpha_0+d^*_p G_p\mathcal{L}_{\phi(t)}^{1,0}\alpha(t)\Big\}.
\end{align*}
The set $\{t\in B\mid\dim V_{t}^{p,k}\geq n\}$ is an analytic subset of $B$ for any $n\in \mathbb{N}$.
\end{definition}

\begin{definition}\label{def-f_t-g_t}
For any $t\in B$, we set
\begin{align*}
f_t: &V_{t}^{p,k} \longrightarrow \frac{F^pA^{k}(X)\cap\ker d_{\phi(t)}\cap\ker d^*_p}{F^pA^{k}(X)\cap\im d_{\phi(t)}\cap\ker d^*_p}\cong F^pH^{k}_{d_{\phi(t)}}(X),\\
&\alpha_0\longmapsto \alpha(t)=\sum_{k\geq0} \alpha_k,~\text{where}~\alpha_{k}=d^*_p G_p\sum_{i+j=k}\mathcal{L}_{\phi_i}^{1,0}\alpha_{j},~\text{for~any}~ k> 0.
\end{align*}
\end{definition}

\begin{proposition}\label{prop V-E-t}
If $V=\ker\triangle_{p}\cap F^{p}A^{k}(X)$, then $f_t$ is surjective with
\[
\ker f_t\cong F^pA^{k}(X)\cap\im d_{\phi(t)}\cap\ker d^*_p.
\]
\end{proposition}
\begin{proof}By Proposition \ref{prop-unique-power-seris-solution}, the following map
\begin{align*}
\tilde{f}_t: &V_{t}^{p,k} \longrightarrow F^pA^{k}(X)\cap\ker d_{\phi(t)}\cap\ker d^*_p,\\
&\alpha_0\longmapsto \alpha(t)=\sum_{k\geq0} \alpha_k,~\text{where}~\alpha_{k}=d^*_p G_p\sum_{i+j=k}\mathcal{L}_{\phi_i}^{1,0}\alpha_{j},~\text{for~any}~ k\neq 0,
\end{align*}
is an isomorphism.
\end{proof}

\begin{theorem}\label{thm-jumping-formula}
For any $0\leq p\leq k$ and $t\in B$, set
\[
v^{p,k}_t:=\dim \ker\triangle_{p}\cap F^{p}A^{k}(X)-\dim F^pA^{k}(X)\cap\ker d_{\phi(t)}\cap\ker d^*_p \geq 0,
\]
then we have $(v^{p,-1}_t:=0)$
\begin{equation}\label{eq-jumping-formula}
\dim F^pH^{k}(X)=\dim F^pH^{k}(X_t)+v^{p,k}_t+v^{p,k-1}_t.
\end{equation}
\end{theorem}
\begin{proof}First, by Proposition \ref{prop V-E-t}, we have
\[
\dim V_{t}^{p,k}-\dim F^pA^{k}(X)\cap\im d_{\phi(t)}\cap\ker d^*_p=\dim F^pH^{k}_{d_{\phi(t)}}(X)=\dim F^pH^{k}(X_t).
\]
Combining this with $(5)$ of Proposition \ref{prop-unique-power-seris-solution}, we have
\[
\dim V_{t}^{p,k}-\dim\ker\triangle_{p}\cap F^{p}A^{k-1}(X)+\dim\left(\ker d_p^*\cap\ker d_{\phi(t)}\right)^{k-1}=\dim F^pH^{k}(X_t),
\]
which by the definition of $v^{p,k-1}_t$ implies
\[
\dim V_{t}^{p,k}=\dim F^pH^{k}(X_t)+v^{p,k-1}_t.
\]
On the other hand, we know that
\begin{equation}\label{eq-Vt-vt}
\dim V_{t}^{p,k}=\dim F^pA^{k}(X)\cap\ker d_{\phi(t)}\cap\ker d^*_p=\dim \ker\triangle_{p}\cap F^{p}A^{k}(X)-v^{p,k}_t.
\end{equation}
Hence, \eqref{eq-jumping-formula} follows.
\end{proof}
Notice that by $(5)$ of Proposition \ref{prop-unique-power-seris-solution}, we have
\begin{equation}\label{eq-vt-imdcapkerd*}
v_{t}^{p,k}=\dim\left( \ker d_p^*\cap\im d_{\phi(t)}\cap F^{p}A^{k+1}(X)\right).
\end{equation}

Forms in $F^pA^{p+q}(X)=\bigoplus_{i=0}^{\infty}A^{p+i,q-i}(X)$ will be called \textit{filtered $(p,q)$-forms on $X$}.
\begin{definition}\label{def-deformation-Dol-Zrpq}
Given
\[
\alpha_0\in\ker d\cap F^{p}A^{p+q}(X),
\]
and $T\subseteq B$, which is an analytic subset of $B$ containing $0$, a \textit{deformation} of $\alpha_0$ (w.r.t. $\pi$ ) on $T$ is a family of forms
\[
\alpha (t)=\alpha^{p,q}(t)+\alpha^{p+1,q-1}(t)+\cdots+\alpha^{n,p+q-n}(t)
\in F^{p}A^{p+q}(X),
\]
such that
\begin{itemize}
  \item[(1)] $\alpha (t)$ is holomorphic in $t$ and $\alpha(0) = \alpha_0$;
  \item[(2)] for any $t\in T$, we have
  \begin{equation}\label{eq-ext-eq}
  d_{\phi(t)}\alpha (t)=0.
  \end{equation}
\end{itemize}
In particular, the Bott-Chern deformation of $(p,q)$-forms can be regarded as deformation of filtered $(p,q)$-forms in this sense, see \cite{Xia19dBC,HX25}.

By a \emph{canonical deformation of $\alpha_0$} w.r.t. $\pi$, we mean a deformation of the form
\begin{equation}\label{eq-canonical}
\alpha(t)=\sum_{k\geq0} \hat{\alpha}_k+d_{\phi(t)}\beta_0,
\end{equation}
where $\beta_0\in F^{p}A^{p+q-1}(X)$ satisfies $d\beta_0=\alpha_0-\mathcal{H}_p\alpha_0$ and $\sum_{k\geq0} \hat{\alpha}_k$ is given by
\[
\hat{\alpha}_0=\mathcal{H}_p\alpha_0,~\hat{\alpha}_k=d^*_p G_p\sum_{i+j=k} \mathcal{L}_{\phi_i}^{1,0}\hat{\alpha}_{j}\in F^{p}A^{p+q}(X),\quad k>0.
\]
We say that the \textit{deformations of filtered $(p,q)$-forms on $X$ are (canonically) unobstructed} (w.r.t. $\pi$) if for any $\alpha_0\in\ker d\cap F^{p}A^{p+q}(X)$
there is a (canonical) deformation of $\alpha_0$ on $B$.

At last, a \emph{deformation $\alpha(t)$ of $\alpha_0$} on $T$ (w.r.t. $\pi$) is said to be \emph{faithful} if $\alpha(t)$ satisfies the following: $\alpha_0\in\im d\cap F^{p}A^{p+q}(X)$ whenever $\alpha(t)\in F^pA^{p+q}(X)\cap\im d_{\phi(t)}$ for some $t\in T$.
\end{definition}
\begin{corollary}\label{Coro-123}
The following holds:
\begin{enumerate}
  \item $\dim F^pH^{k}(X_t)$ is independent of $t\in B$ if and only if the deformations of filtered $(p,k-p)$-forms and filtered $(p,k-1-p)$-forms are canonically unobstructed;
  \item For any $0\leq p\leq k$, the alternating sum $\sum_{i=0}^k (-1)^{k-i}\dim F^pH^{i}(X_t)$ and $\dim E^{0,q}_1(X_t)/E^{0,q}_\infty(X_t)$ are upper semi-continuous functions of $t\in B$ (in analytic Zariski topology). Moreover, the deformations of filtered $(p,k-p)$-forms are canonically unobstructed if and only if $\sum_{i=0}^k (-1)^{k-i}\dim F^pH^{i}(X_t)$ is independent of $t$;
  \item The notion of canonically unobstructedness for filtered forms are independent of the choices of Hermitian metrics.
\end{enumerate}
\end{corollary}
\begin{proof}$(1)$ By the definition of canonically unobstructedness, the deformations of filtered $(p,k-p)$-forms are canonically unobstructed iff
\[
V_{t}^{p,k}=\ker\triangle_{p}\cap F^{p}A^{k}(X),\quad\text{for~any}~t\in B,
\]
which is equivalent to $v^{p,k}_t=0$ for any $t\in B$. The conclusion then follows from \eqref{eq-jumping-formula}.

$(2)$ By considering the alternating sum of \eqref{eq-jumping-formula}, we have
\[
v^{p,k}_t=\sum_{i=0}^k \dim (-1)^{k-i}F^pH^{i}(X)-\sum_{i=0}^k \dim (-1)^{k-i}F^pH^{i}(X_t),
\]
which combined with \eqref{eq-Vt-vt} implies
\[
\sum_{i=0}^k \dim (-1)^{k-i}F^pH^{i}(X_t)=\dim V_{t}^{p,k}.
\]
The upper semi-continuity of $\sum_{i=0}^k (-1)^{k-i}\dim F^pH^{i}(X_t)$ follows from this and the fact that $\dim V_{t}^{p,k}$ is upper semi-continuous. As a consequence, the upper semi-continuity of $\dim E^{0,q}_1(X_t)/E^{0,q}_\infty(X_t)$ follows from \eqref{eq-E0q}.

$(3)$ follows directly from $(2)$ since $\dim F^pH^{i}(X_t)$ are complex structure invariants.
\end{proof}

\begin{proposition}\label{prop-faithful}
Assume the deformations of filtered $(p,k-1-p)$-forms are canonically unobstructed. Let $\alpha(t)$ be a canonical deformation of $\alpha_0\in\ker d\cap F^{p}A^{k}(X)$ on $T$ (w.r.t. $\pi$), then it is faithful.
\end{proposition}
\begin{proof}
We may assume $\alpha(t)$ is given by \eqref{eq-canonical}. If $\alpha(t)\in F^pA^{k}(X)\cap\im d_{\phi(t)}$ for some $t\in T$, then we have $\sum_{m\geq0} \hat{\alpha}_m\in \im d_{\phi(t)}\cap\ker d_p^*$. But since the deformations of filtered $(p,k-1-p)$-forms are canonically unobstructed, we have $v^{p,k-1}_t=0$ for any $t\in B$ which implies (in view of \eqref{eq-vt-imdcapkerd*})
\[
\ker d_p^*\cap\im d_{\phi(t)}\cap F^{p}A^{k}(X)=0,\quad\text{for~any}~t\in B.
\]
Hence, we get $\sum_{m\geq0} \hat{\alpha}_m=0\Rightarrow\hat{\alpha}_0=\mathcal{H}_p\alpha_0=0$ by using $(3)$ of Proposition \ref{prop-unique-power-seris-solution}.
\end{proof}

\section{Deformations of Hodge filtrations}
In this section, we prove the main result of this paper which concerns geometric properties of period maps for possibly non-K\"ahler manifolds.
\subsection{Faithful deformations of filtered forms}
First of all, we consider the problem of when the deformations of filtered $(p,q)$-forms are canonically unobstructed. In fact, we will use the following result whose proof will be sketched below for the reader's convenience.
\begin{theorem}{\cite[Lem.\,3.3]{WX23}}\label{thm-WX23} Let $\pi: (\mathcal{X}, X)\to (B,0)$ be a complex analytic family over a small polydisc $B\subset\C^m$ such that for each $t\in B$ the complex structure on $X_t$ is represented by a Beltrami differential $\phi(t)$. Assume
  \[
  F^{p}A^{p+q+1}(X)\cap dA^{p+q}(X) =dF^{p}A^{p+q}(X).
  \]
Then the deformations of filtered $(p,q)$-forms are canonically unobstructed.
\end{theorem}
\begin{proof}For any given $\alpha_0\in F^{p}A^{p+q}(X)\cap\ker d$, we first try to a (formal) power series solution
\[
\alpha(t)=\sum_{k\geq0}\alpha_k,\quad \text{with}~\alpha_k\in F^{p}A^{p+q}(X)~\text{for~any}~k,
\]
to the following
\[
d_{\phi(t)}\alpha(t)=0,
\]
or
\begin{equation}\label{eq-induction}
d\alpha_k=\sum_{1\leq j\leq k}\mathcal{L}_{\phi_j}^{1,0}\alpha_{k-j},\quad   k\geq 0,
\end{equation}
where for $k=0$, we set $\sum_{1\leq j\leq k}\mathcal{L}_{\phi_j}^{1,0}\alpha_{k-j}=0$. Assume this can be solved for $1\leq k\leq N$, it follows from \cite[Lem.\,3.2]{WX23} that
\begin{equation}
\left(d(e^{i_{\phi(t)}}-1)\alpha(t)\right)_{N_1}=-\left(\mathcal{L}_{\phi(t)}^{1,0}\alpha(t)\right)_{N_1},\quad \forall N_1\leq N+1.
\end{equation}
This combined with our assumption means that
\[
\sum_{j=1}^{N+1} \mathcal{L}_{\phi_j}^{1,0} \alpha_{N+1-j}\in F^{p}A^{p+q+1}(X)\cap dA^{p+q}(X) =dF^{p}A^{p+q}(X).
\]
Hence, there is a $\alpha_{N+1}\in F^{p}A^{p+q}(X)$ such that
\[
  d \alpha_{N+1}=\sum_{j=1}^{N+1} \mathcal{L}_{\phi_j}^{1,0} \alpha_{N+1-j}.
\]
We have proved that \eqref{eq-induction} can be solved inductively. Moreover, in each step of the induction, we could have chosen the canonical solution to equation \eqref{eq-induction}, namely,
\[
\alpha_k:=d^*_p G_p\sum_{i+j=k} \mathcal{L}_{\phi_i}^{1,0}\alpha_{j}\in F^{p}A^{p+q}(X),\quad k>0.
\]
Because $\alpha_0\in F^{p}A^{p+q}(X)\cap\ker d$ is arbitrary, this shows the deformations of filtered $(p,q)$-forms are canonically unobstructed.
\end{proof}

As a consequence, we have
\begin{theorem}\label{thm-exp-FpH}
Let $\pi: (\mathcal{X}, X)\to (B,0)$ be a complex analytic family over a small polydisc $B\subset\C^m$ such that for each $t\in B$ the complex structure on $X_t$ is represented by a Beltrami differential $\phi(t)$. Fix $p\leq k$ and assume
  \[
  F^{p}A^{k+1}(X)\cap dA^{k}(X) =dF^{p}A^{k}(X),
  \]
and
  \[
  F^{p}A^{k}(X)\cap dA^{k-1}(X) =dF^{p}A^{k-1}(X).
  \]
For any $t\in B$, the exponential operator induces the following isomorphism of vector spaces (still denoted by $e^{i_{\phi(t)}}$)
\begin{equation}\label{eq-exponential-coho}
e^{i_{\phi(t)}}: F^pH^{k}(X)\longrightarrow F^pH^{k}(X_t):[\alpha_0]\longmapsto [e^{i_{\phi(t)}}\alpha(t)],
\end{equation}
where $\{\alpha(t)\}_{t\in B}$ is the canonical deformation of $\alpha_0$.
\end{theorem}
\begin{proof}First of all, note that \eqref{eq-exponential-coho} is the composition of the following two mappings
\begin{equation}\label{eq-can-def-map}
F^{p}H^{k}(X)\longrightarrow F^{p}H^{k}_{d_{\phi(t)}}(X):[\alpha_0]\longmapsto [\alpha(t)],
\end{equation}
where $\alpha(t)$ is the canonical deformation of $\alpha_0$ and
\begin{equation}\label{eq-dphi-dXt}
F^{p}H^{k}_{d_{\phi(t)}}(X)\longrightarrow F^{p}H^{k}(X_t):[\alpha]\longmapsto [e^{i_{\phi(t)}}\alpha].
\end{equation}
First, according to Theorem \ref{thm-WX23}, the deformations of filtered $(p,k-p)$-forms and filtered $(p,k-1-p)$-forms are canonically unobstructed. By Proposition \ref{prop-faithful}, all canonical deformations of filtered $(p,k-p)$-forms in this case are faithful. These two facts imply the homomorphism given by \eqref{eq-can-def-map} is well-defined and injective. The homomorphism given by \eqref{eq-dphi-dXt} is an isomorphism in view of \eqref{eq-wdc-dc}. So we see that \eqref{eq-exponential-coho} is well-defined and injective. In particular, we have
\[
\dim F^pH^{k}(X)\leq \dim F^pH^{k}(X_t),\quad\text{for~any}~t\in B.
\]
But according to Proposition \ref{prop-upper-semicontinuous}, $\dim F^pH^{\bullet}(X_t)$ is upper semi-continuous which implies
\[
\dim F^pH^{k}(X)= \dim F^pH^{k}(X_t),\quad\text{for~any}~t\in B.
\]
Therefore, \eqref{eq-exponential-coho} is an isomorphism.
\end{proof}

\subsection{Geometry of the period mappings}
In the context of Theorem \ref{thm-exp-FpH}, the \emph{period map} can be defined as follows
\begin{equation}\label{eq-period-map}
\Phi^{p,k}:B\longrightarrow \G (f^{p,k}, H^{k}(X,\C)),\quad t\longmapsto F^pH^{k}(X_t),
\end{equation}
where $f^{p,k}:=\dim F^pH^{k}(X)$. The period map $\Phi^{p,k}$ has the following properties (cf.~\cite[Chp.\,10]{Voi02I} for the case of K\"ahler manifolds),

\begin{theorem}\label{thm-period-map}
Fix $p\leq k$ and assume
  \[
  F^{p}A^{k+1}(X)\cap dA^{k}(X) =dF^{p}A^{k}(X),
  \]
and
  \[
  F^{p}A^{k}(X)\cap dA^{k-1}(X) =dF^{p}A^{k-1}(X).
  \]
The following holds:
\begin{enumerate}
  \item The period map $\Phi^{p,k}$ is holomorphic;
  \item Griffiths transversality: the tangent map
  \[
  d\Phi^{p,k}_0:T_0B\longrightarrow  \Hom(F^pH^{k}(X), H^{k}(X,\C)/F^pH^{k}(X))
  \]
  has values in $\Hom(F^pH^{k}(X), F^{p-1}H^{k}(X)/F^pH^{k}(X))$;
\end{enumerate}
Let now $\omega$ be a Hermitian metric on $X$. Denote by $\mathcal{H}_p$ the orthogonal projection onto the space of harmonic filtered forms $\ker\triangle_p$ in $F^pA^\bullet(X)$ (see Subsection \ref{subsec-Hodge-decomp}). If $\mathcal{H}_p=\mathcal{H}_{p-1}\mid_{F^{p}A^{\bullet}(X)}$ and $\dim B=1$, then the tangent map of $\Phi^{p,k}$ is given by
\[
d\Phi^{p,k}_0(\frac{\p}{\p t}\Big|_{t=0})([\alpha_0])=[\mathcal{H}_{p-1}i_{\phi_1}\alpha_0]\in F^{p-1}H^{k}(X),\quad [\alpha_0]\in F^pH^{k}(X),
\]
where $t$ is a local coordinate function on $B$ and $\phi_1$ is the first order coefficient in $\phi(t)=\sum_{j\geq 1}\phi_jt^j$.
\end{theorem}
\begin{proof}
These statements can be deduced from Theorem \ref{thm-exp-FpH}. In fact, for any $t\in B$ we now have
\[
\Phi^{p,k}(t)=F^pH^{k}(X_t)=\C\Big\{ [e^{i_{\phi(t)}}\alpha^l(t)] \Big\}_{l=1}^{f^{p,k}},
\]
where $\{[\alpha^l_0]: l=1,\cdots,f^{p,k}\}$ is a basis of $F^pH^{k}(X)$
and $\alpha^l(t)$ is the canonical deformation of $\alpha^l_0$. This shows $\Phi^{p,k}$ is holomorphic because both $\phi(t)$ and each $\alpha^l(t)$ are holomorphic. This is $(1)$.

Next, let $\alpha(t)$ be the canonical deformation of some $[\alpha_0]\in F^pH^{k}(X)$ such that $\alpha_0\in\ker\triangle_{p}\cap F^{p}A^{k}(X)$, i.e. $\alpha_0$ is a harmonic filtered $(p,k-p)$-form. Without loss of generality assume $\dim B=1$ and we can write $\phi(t)=\sum_{j\geq 1}\phi_jt^j$, $\alpha(t)=\sum_{m\geq 0}\alpha_mt^m$, where $\alpha_{m}=d^*_p G_p\sum_{i+j=m}\mathcal{L}_{\phi_i}^{1,0}\alpha_{j}$. So we have
\[
\frac{\p}{\p t}\Big|_{t=0}\alpha(t)=\alpha_1=d^*_p G_p\mathcal{L}_{\phi_1}^{1,0}\alpha_{0}~.
\]
Now if we differentiate $e^{i_{\phi(t)}}\alpha(t)$ with respect to $t$ at $t=0$, we get
\begin{equation}\label{eq-differentiate}
\frac{\p}{\p t}\Big|_{t=0}[e^{i_{\phi(t)}}\alpha(t)]=[i_{\phi_1}\alpha_0+d^*_p G_p\mathcal{L}_{\phi_1}^{1,0}\alpha_{0}]\in F^{p-1}H^{k}(X).
\end{equation}
This is $(2)$.

Finally, if $\dim B=1$ then we can write $\phi(t)=\sum_{j\geq 1}\phi_jt^j$ with $t\in B\subset\C$ and
\[
\kappa(\frac{\p}{\p t}\Big|_{t=0})=\frac{\p}{\p t}\Big|_{t=0}\phi(t)=\phi_1.
\]
Let $[\alpha_0]\in F^pH^{k}(X)$ and $\alpha(t)=\sum_{m\geq 0}\alpha_mt^m$ a canonical deformation of $\alpha_0\in F^pA^{k}(X)\cap\ker d$. From $(1),(2)$, we know that locally,
\begin{align*}
d\Phi^{p,k}_0(\frac{\p}{\p t}\Big|_{t=0})([\alpha_0])
&=\frac{\p}{\p t}\Big|_{t=0}[e^{i_{\phi(t)}}\alpha(t)]\\
&=[i_{\phi_1}\alpha_0+d^*_p G_p\mathcal{L}_{\phi_1}^{1,0}\alpha_{0}]\\
&=[\mathcal{H}_{p-1}i_{\phi_1}\alpha_0+\mathcal{H}_{p-1}d^*_p G_p\mathcal{L}_{\phi_1}^{1,0}\alpha_{0}]\\
&=[\mathcal{H}_{p-1}i_{\phi_1}\alpha_0+\mathcal{H}_{p}d^*_p G_p\mathcal{L}_{\phi_1}^{1,0}\alpha_{0}]\\
&=[\mathcal{H}_{p-1}i_{\phi_1}\alpha_0]\in F^{p-1}H^{k}(X).
\end{align*}
\end{proof}

Note that $\{e^{i_{\phi(t)}}\alpha^l(t)\}_{l=1}^{f^{p,k}}$ is a holomorphic frame of the Hodge bundle $\mathcal{F}^p:=\bigcup_{t\in B}F^pH^{k}(X_t)$ (cf. \cite{LS18}).

\subsection{The case of K\"ahler manifolds}
Finally, we make some remarks on the situation of K\"ahler manifolds.

If $X$ is equipped with a K\"ahler metric, the period map
\[
\Phi^{p,k}:B\longrightarrow \G (f^{p,k}, H^{k}(X,\C)),
\]
defined in \eqref{eq-period-map} can be alternatively described by using canonical Dolbeault deformations of $(p,k-p)$-forms (cf. \cite{Xia22}) instead of using canonical deformations of filtered forms as in the proof of Theorem \ref{thm-period-map}. In fact, let $\{[\alpha^l_0]: l=1,\cdots,f^{p,k}\}$ be a basis of $F^pH^{k}(X)$ such that
\begin{align*}
\mathcal{H}^{p,k-p}(X)=&\C\{\alpha^1_0,\cdots, \alpha^{h^{p,k-p}}_0\},\\
\mathcal{H}^{p+1,k-p-1}(X)=&\C\{\alpha^{h^{p,k-p}+1}_0,\cdots, \alpha^{h^{p,k-p}+h^{p+1,k-p-1}}_0\},\cdots,
\end{align*}
where $\mathcal{H}^{p,k-p}(X)$ is the space of $\pb$-harmonic $(p,k-p)$-forms and $h^{p,k-p}=\dim\mathcal{H}^{p,k-p}(X)$ is the $(p,k-p)$ Hodge number.
Then locally the period map $\Phi^{p,k}$ is given by
\[
\Phi^{p,k}(t)=F^pH^{k}(X_t)=\C\Big\{ [e^{i_{\phi(t)}}\alpha^l(t)] \Big\}_{l=1}^{f^{p,k}},
\]
where for each $i\geq0$, $\alpha^{p+i,k-p-i}(t):=(\alpha^l(t))^{p+i,k-p-i}$ is the canonical Dolbeault deformation of $(\alpha^l_0)^{p+i,k-p-i}$. Here $(\bullet)^{p+i,k-p-i}$ denotes the $(p+i,k-p-i)$-th component of the $k$-form $\bullet$. More precisely, each $\alpha^{p+i,k-p-i}(t)=\sum_{m\geq0}\alpha_m^{p+i,k-p-i}t^m$ is a convergent power series with
\begin{itemize}
  \item $\alpha_0^{p+i,k-p-i}=(\alpha^l_0)^{p+i,k-p-i}$;
  \item $\alpha_m^{p+i,k-p-i}=\sum_{j\geq1}\pb^*G_{\pb}
      \mathcal{L}_{\phi_j}^{1,0}\alpha_{m-j}^{p+i,k-p-i}
      =\footnote{This follows from the K\"ahler identity.}-\p\sum_{j\geq1}\pb^*G_{\pb}
      i_{\phi_j}\alpha_{m-j}^{p+i,k-p-i}$;
  \item $\alpha^{p+i,k-p-i}(t)\in\ker\pb_{\phi(t)}\cap\ker\p$ which implies $e^{i_{\phi(t)}}\alpha^l(t)\in\ker d$ in view of \eqref{eq-LRY15}.
\end{itemize}
So \eqref{eq-differentiate} now becomes
\begin{equation}\label{eq-differentiate-kahler}
\frac{\p}{\p t}\Big|_{t=0}[e^{i_{\phi(t)}}\alpha(t)]=[i_{\phi_1}\alpha_0+\alpha_1]=
[i_{\phi_1}\alpha_0-\p\pb^*G_{\pb}i_{\phi_1}\alpha_0]=[i_{\phi_1}\alpha_0],
\end{equation}
where in the last equality we have used the fact that for K\"ahler manifolds, the $d,\p,\pb$-harmonic spaces all coincide so that the harmonic projection operators $\mathcal{H}_d, \mathcal{H}_\p, \mathcal{H}_\pb$ are the same which implies
\[
\mathcal{H}_d(i_{\phi_1}\alpha_0-\p\pb^*G_{\pb}i_{\phi_1}\alpha_0)
=\mathcal{H}_di_{\phi_1}\alpha_0-\mathcal{H}_\p\p\pb^*G_{\pb}i_{\phi_1}\alpha_0
=\mathcal{H}_di_{\phi_1}\alpha_0=\mathcal{H}_\pb i_{\phi_1}\alpha_0.
\]
In particular, it follows from \eqref{eq-differentiate-kahler} that the following diagram is commutative (see e.g. \cite[Thm.\,10.21]{Voi02I}):
 \[
  \xymatrix{
  T_0B \ar[r]^-{\kappa}\ar[dr]_-{d\Phi^{p,k}_0} &  H^{1}(X,T_X^{1,0})  \ar[d]^-{\iota}   \\
          & \Hom(F^pH^{k}(X), F^{p-1}H^{k}(X)/F^pH^{k}(X)) , }
 \]
where $\kappa$ is the Kodaira-Spencer map and $\iota$ is defined as follows: for any $[\varphi]\in H^{1}(X,T_X^{1,0})$, $\iota([\varphi])=i_\varphi$. In the situation of Theorem \ref{thm-period-map}, the homomorphism $\iota$ may not be well-defined.

To end this section, we next consider the condition $\mathcal{H}_p=\mathcal{H}_{p-1}\mid_{F^{p}A^{\bullet}(X)}$ which appeared in Theorem \ref{thm-period-map}. We will show that it automatically holds for K\"ahler manifolds.

For this, we first make some observations about the operator $\triangle_p$. In Subsection \ref{subsec-Hodge-decomp}, we have introduced the Laplacian operator for filtered forms as the following
\[
\triangle_p=dd^*_p+d^*_pd:F^{p}A^{\bullet}(X)\longrightarrow F^{p}A^{\bullet}(X),
\]
where $d^*_p=\Pi^{\geq p}d^*$. The operator $\triangle_p$ is related to the usual Laplacian $\triangle=dd^*+d^*d$ in the following way:
\begin{lemma}\label{lem} For any $\alpha=\alpha^{p,q}+\alpha^{p+1,q-1}+\cdots\in F^{p}A^{p+q}(X)$, we have
  \begin{equation}\label{eq-Lapp-Lap}
   \triangle_p\alpha = \triangle\alpha -[\pb,\p^*]\alpha^{p,q}-\p\p^*\alpha^{p,q},
  \end{equation}
where $[\pb,\p^*]=\pb\p^*+\p^*\pb$.
Furthermore, if the Hermitian metric on $X$ is K\"ahler, then
  \begin{equation}\label{eq-Lapp-Lap-Kahler}
  \triangle_p\alpha = \triangle\alpha -\p\p^*\alpha^{p,q}.
  \end{equation}
\end{lemma}
\begin{proof}First of all, if $\alpha\in F^{p+1}A^{p+q}(X)$, then it is clear that
\[
  \triangle_p\alpha = \triangle\alpha.
\]
It follows that for any $\alpha\in F^{p}A^{p+q}(X)$,
\begin{align*}
\triangle_p\alpha
&= \triangle_p(\alpha^{p,q}+\alpha-\alpha^{p,q})\\
&=\triangle_p\alpha^{p,q}+\triangle(\alpha-\alpha^{p,q})\\
&=(dd^*_p+d^*_pd)\alpha^{p,q}+\triangle(\alpha-\alpha^{p,q})\\
&=d\pb^*\alpha^{p,q}+d^*d\alpha^{p,q}-\p^*\pb\alpha^{p,q}+\triangle(\alpha-\alpha^{p,q})\\
&=\triangle\alpha^{p,q}-d\p^*\alpha^{p,q}-\p^*\pb\alpha^{p,q}
+\triangle(\alpha-\alpha^{p,q})\\
&=\triangle\alpha -[\pb,\p^*]\alpha^{p,q}-\p\p^*\alpha^{p,q}.
\end{align*}
If the Hermitian metric on $X$ is K\"ahler, then $[\pb,\p^*]=0$ (see e.g. \cite[Chp.\,VI]{Dem12}).
\end{proof}

\begin{proposition}\label{prop-kahler} If the Hermitian metric on $X$ is K\"ahler, then
  \begin{equation}\label{eq-Hp=Hp-1}
  \mathcal{H}_{p+1}=\mathcal{H}_{p}\mid_{F^{p+1}A^{\bullet}(X)},\quad p\geq0.
  \end{equation}
\end{proposition}
\begin{proof}
We first show $\ker\triangle_p\cap F^{p+1}A^{\bullet}(X)=\ker\triangle_{p+1}$ holds for any $p\geq0$. Indeed, for any $\alpha\in \ker\triangle_p\cap F^{p+1}A^{\bullet}(X)$, we have $d\alpha=0$ and $d^*\alpha=d_p^*\alpha=0$ which implies $d_{p+1}^*\alpha=0$. Hence, $\ker\triangle_p\cap F^{p+1}A^{\bullet}(X)\subseteq\ker\triangle_{p+1}$. On the other hand, working in total degree $p+q$, let $\alpha\in \ker\triangle_{p+1}\cap F^{p+1}A^{p+q}(X)$. By using \eqref{eq-Lapp-Lap-Kahler} in Lemma \ref{lem} we have (since $\triangle=2\triangle_\p=2\triangle_\pb$)
\begin{align*}
0=\triangle_{p+1}\alpha =& \triangle\alpha -\p\p^*\alpha^{p+1,q-1}\\
=&(2\triangle_\p\alpha^{p+1,q-1}-\p\p^*\alpha^{p+1,q-1})+\triangle(\alpha-\alpha^{p+1,q-1})\\
\in& A^{p+1,q-1}(X)\oplus F^{p+2}A^{\bullet}(X),
\end{align*}
which implies
\[
 \p\p^*\alpha^{p+1,q-1}=0\Longrightarrow\triangle\alpha=0
 \Longrightarrow \triangle_p\alpha=\triangle\alpha=0,
\]
where $\triangle_p\alpha=\triangle\alpha$ follows from \eqref{eq-Lapp-Lap-Kahler} and the fact that $\alpha^{p,q}=0$. Therefore, we get $\ker\triangle_{p+1}\subseteq\ker\triangle_p\cap F^{p+1}A^{\bullet}(X)$.

The preceding argument also shows that
\begin{equation}\label{eq-laplacian-r}
\ker\triangle_r=\ker\triangle\cap F^rA^\bullet(X),
\qquad r\geq0.
\end{equation}
Since on a K\"ahler manifold the ordinary harmonic projection $\mathcal H=\mathcal H_0$ preserves bidegrees,
we have $\mathcal H\big(F^rA^\bullet(X)\big)\subset F^rA^\bullet(X)$. So $\mathcal H_r=\mathcal H|_{F^rA^\bullet(X)}$. Indeed, for any $\alpha\in F^rA^\bullet(X)$, we have $\triangle G\alpha\perp\ker\triangle_r$ in view of \eqref{eq-laplacian-r}. On the other hand, both $\mathcal H\alpha$ and $\triangle G\alpha$ belong to $F^rA^\bullet(X)$. It follows that $\alpha=\mathcal H\alpha+\triangle G\alpha$ is also the Hodge decomposition of $\alpha$ with respect to $\triangle_r$.

The asserted identity \eqref{eq-Hp=Hp-1} now follows.
\end{proof}

\vskip 1\baselineskip \textbf{Declaration of generative AI.}
During the preparation of this work, the author used Codex to assist with proofreading and manuscript revision. The author reviewed and edited the output as needed and takes full responsibility for the content of the published article.

\bibliographystyle{alpha}
\bibliography{reference}
\end{document}